\documentclass[9pt,twocolumn,twoside]{pnas-new}

\templatetype{pnasresearcharticle} 

\newtheorem{theorem}{Theorem}
\newtheorem{assumption}{Assumption}
\newtheorem{lemma}{Lemma}

\newtheorem{definition}{Definition}

\begin{document}

\title{Quantization Avoids Saddle Points in Distributed Optimization }

\author[a]{Yanan Bo}
\author[a,1]{Yongqiang Wang}

\affil[a]{Department of Electrical and Computer Engineering, Clemson University, Clemson, SC 29634 USA}


\leadauthor{Bo}

\significancestatement{
Distributed optimization underpins key functionalities of numerous engineered systems such as smart grids, intelligent transportation, and smart cities. It is also reshaping the landscape of machine learning due to its inherent advantages in handling large data/model sizes. However, saddle-point avoidance becomes extremely challenging in distributed optimization because individual agents in distributed optimization do not have access to the global gradient. We show that quantization effects, which are unavoidable due to communications in distributed optimization and regarded as detrimental in existing studies, can be exploited to enable saddle-point avoidance for free. By judiciously designing the quantization scheme, we propose an approach that evades saddle points and ensures convergence to a second-order stationary point in distributed nonconvex optimization.
}

\authorcontributions{Y.B and Y.W designed research, performed research and wrote the paper.}
\authordeclaration{The authors declare no competing interest.}
\correspondingauthor{\textsuperscript{1}To whom correspondence should be addressed. E-mail: yongqiw@clemson.edu }


\begin{abstract}
Distributed nonconvex optimization underpins key functionalities of numerous distributed systems, ranging from power systems, smart buildings, cooperative robots, vehicle networks to sensor networks. Recently, it has also merged as a promising solution to handle the enormous growth in data and model sizes in deep learning. A fundamental problem in distributed nonconvex optimization is avoiding convergence to saddle points, which significantly degrade optimization accuracy. We discover that the process of quantization, which is necessary for all digital communications, can be exploited to enable saddle-point avoidance. More specifically, we propose a stochastic quantization scheme and prove that it can effectively escape saddle points and ensure convergence to a second-order stationary point in distributed nonconvex optimization. With an easily adjustable quantization granularity, the approach allows a user to control the number of bits sent per iteration and, hence, to aggressively reduce the communication overhead. Numerical experimental results using distributed optimization and learning problems on benchmark datasets confirm the effectiveness of the approach.
\end{abstract}

\dates{This manuscript was compiled on \today}
\doi{\url{www.pnas.org/cgi/doi/10.1073/pnas.XXXXXXXXXX}}

\maketitle
\thispagestyle{firststyle}
\ifthenelse{\boolean{shortarticle}}{\ifthenelse{\boolean{singlecolumn}}{\abscontentformatted}{\abscontent}}{}

\firstpage[3]{5}




With the unprecedented advances in embedded electronics and communication technologies, cooperation or coordination has emerged as a key feature in numerous engineered systems such as smart grids, intelligent transportation systems, cooperative robots, cloud computing, and smart cities. This has spurred the development of distributed algorithms in which spatially distributed computing devices (hereafter referred to as agents), communicating over a network, cooperatively solve a task without resorting to a central coordinator/mediator that aggregates all data in the network. Many of these distributed algorithms boil down to the following distributed optimization problem:

\begin{equation}\label{Equ:P1}
\begin{aligned}
\min_{\boldsymbol{\theta}\in \mathbb{R}^{d}} F(\boldsymbol{\theta})=\frac{1}{N}\sum_{i=1}^{N}f_i(\boldsymbol{\theta}),
\end{aligned}
\end{equation}
where $ f_i(\cdot): \mathbb{R}^d \rightarrow \mathbb{R} $ denotes the local objective function private to agent $i$, $ F(\cdot): \mathbb{R}^d \rightarrow \mathbb{R} $ is the global objective function representing the network-level cost to be minimized cooperatively by all participating agents, and $N$ is the number of agents.

Initially introduced in the 1980s in the context of parallel and distributed computation \cite{bertsekas2015parallel}, the above distributed optimization problem has received intensive interest in the past decade due to the surge of smart systems and deep learning applications \cite{nedic2009distributed, nedic2010constrained}. So far, plenty of approaches have been proposed to solve the above
distributed optimization problem, with some of the commonly used approaches including gradient methods (see, e.g., \cite{nedic2009distributed,srivastava2011distributed,shi2015extra,xu2017convergence,qu2017harnessing,xin2018linear}), distributed alternating direction method of multipliers (see, e.g., \cite{shi2014linear}), and distributed Newton methods (see, e.g., \cite{zhang2018admm}).

However, most of these approaches
focus on convex objective functions, whereas results are relatively sparse for nonconvex objective functions. In fact, in many applications, the objective functions are essentially nonconvex. For example, in the resource allocation problem of communication networks, the utility functions are nonconvex when the communication traffic is non-elastic \cite{tychogiorgos2013non}; in most machine learning applications, the objective functions are nonconvex due to the presence of multi-layer neural networks \cite{tsianos2012consensus}; in policy optimization for linear-quadratic regulators \cite{fazel2018global} as well as for robust and risk-sensitive control \cite{doi:10.1137/20M1347942}, nonconvex optimization naturally arises.

In nonconvex optimization, oftentimes, the most fundamental problem is to avoid saddle points (stationary points that are not local extrema). For example, in machine learning applications, it has been shown that the main bottleneck in parameter optimization is not due to the existence of multiple local minima but the existence of many saddle points that trap gradient updates \cite{ge2015escaping}. The problem of saddle points is more acute in deep neural networks, where saddle points are usually encircled by high-error plateaus, exerting substantial deceleration on the learning process while engendering a deceptive semblance of the presence of a local minimum \cite{dauphin2014identifying, choromanska2015loss}. To escape saddle points, classical approaches resort to second-order information, in particular, the Hessian matrix of second-order derivatives (see, e.g., \cite{nesterov2006cubic,curtis2017trust}). The Hessian matrix-based approach, however, incurs high costs in both computation and storage. This is because the dimension of the Hessian matrix increases quadratically with an increase in the optimization-variable dimension, which can scale to hundreds of millions in modern deep learning applications \cite{tang2020communication}. Recently,  random perturbations of first-order gradient methods have been shown capable of escaping saddle points in centralized optimization (see, e.g., \cite{ge2015escaping,du2017gradient}). However, it is unclear if this is still true in decentralized nonconvex optimization, where the decentralized architecture brings in fundamental differences in optimization dynamics. For example, in decentralized optimization, the saddle points of individual local objective functions $f_i(\cdot)$ are different from those of the global
objective function $F(\cdot)$, which is the only function that needs to be considered in centralized optimization. In fact, in distributed optimization, all local objective functions $f_i(\cdot)$ are private to individual agents, preventing any single agent from accessing the global objective function $F(\cdot)$ and further from exploiting the gradient/Hessian information of $F(\cdot)$ in its local iteration to avoid the saddle points of $F(\cdot)$.
In addition, the inter-agent coupling also complicates the optimization dynamics. Note that random algorithm initialization has been shown to be able to asymptotically avoid saddle points in centralized nonconvex optimization \cite{lee2016gradient}, which has been further extended to the decentralized case in \cite{daneshmand2020second}. However, the result in \cite{du2017gradient} shows that this approach to avoiding saddle points may take an exponentially longer time, rendering it impractical.

In this paper, we propose to exploit the effects of quantization, which are naturally inherent to all digital communication methods, to evade saddle points in distributed nonconvex optimization. The process of quantization is necessary in all modern communications to represent continuous-valued variables with a smaller set of discrete-valued variables since digital communication channels can only transmit/receive bit streams. The conversion from continuous-valued variables to discrete-valued variables inevitably leads to rounding and truncation errors. In fact, in distributed learning for deep neural networks, since model parameters or gradients have to be shared across agents in every iteration and the dimension of these model parameters and gradients can easily scale to hundreds of millions \cite{tang2020communication}, it is a common practice to use coarse quantization schemes or compression techniques to reduce the overhead of communication \cite{wen2017terngrad, alistarh2017qsgd}.  Recently, plenty of distributed optimization and learning algorithms have been proposed that can ensure provable convergence to the optimal solution in the convex case (see, e.g., \cite{Kashyap2006QuantizedC, Rabbat2005QuantizedIA,el2016design, wang2011control,zhu2011distributed,kia2015distributed, su2016fault,jia2015fault,wang2020adaptive}) or to first-order stationary points in the nonconvex case (see, e.g., \cite{koloskova2019choco,yuan2022revisiting,8320863}), even in the presence quantization/compression errors. However, in all these existing results, quantization effects are treated as detrimental to the distributed optimization process and have to be suppressed to ensure convergence accuracy.  In this paper, to the contrary, we exploit quantization effects to evade saddle points and hence improve convergence accuracy in distributed nonconvex optimization. By judiciously designing the quantization scheme, we propose an algorithm that can make use of quantization effects to effectively escape saddle points and ensure convergence to second-order stationary points. To the best of our knowledge, this is the first time that quantization is shown to be beneficial to the convergence accuracy of distributed optimization. The proposed quantization scheme can also aggressively reduce the overhead of communication, which is widely regarded as the bottleneck in distributed training of machine-learning models \cite{wen2017terngrad}.
\section*{Problem Formulation}\label{Sec:problem}

\subsection*{Notations}
We use bold letters to denote matrices and vectors, i.e., $\boldsymbol{A}$ and $\boldsymbol{x}$.
We use $\|\cdot\|$ to represent the $\ell_2$ norm of vectors and the Frobenius norm of matrices. For a function $ F(\cdot): \mathbb{R}^d \rightarrow \mathbb{R} $, we use $\nabla F (\cdot)$ and $\nabla^2 F(\cdot) $ to denote its gradient and Hessian, respectively. We use $\mathcal{O}(\cdot)$ to hide absolute constants that do not depend on any problem parameter. We use $[N]$ to represent the set $\{1,2,\cdots,N\}$. We use $\lambda_{\min}(\cdot)$ to represent the minimal eigenvalue of a matrix.
\subsection*{Formulation}
We consider a distributed optimization problem where $N$ agents, each with its own local objective function, collaboratively optimize the network-level sum (average) of all local objective functions. Since the local objective functions are private to individual agents, no agents have access to the global objective function. To solve the distributed optimization problem, individual agents have to share local intermediate optimization variables with their respective immediate neighboring agents to ensure convergence to a desired solution.  We describe the local interaction among agents using a weight matrix $\boldsymbol{A}=[a_{ij}]_{N \times N}$, where $a_{ij}>0$ if agent $j$ and agent
$i$ can directly communicate with each other, and $a_{ij} = 0$ otherwise. For an agent $ i \in  [N]$, its neighbor set $\mathcal{N}_i$ is
defined as the collection of agents $j$ such that  $a_{ij}>0$. $a_{ii}$ represents self-interaction, i.e., the influence of agent $i$'s optimization variable at iteration $k$ on its optimization variable at iteration $k+1$. Furthermore, we make the following assumption on $\boldsymbol{A}$:

\begin{assumption}\label{Assum:matrix_A}
    The matrix $\boldsymbol{A}=\{a_{ij}\} \in \mathbb{R}^{N \times N}$ is symmetric and satisfies $ \boldsymbol{1}^\top \boldsymbol{A} = \boldsymbol{1}^\top$, $\boldsymbol{A}\boldsymbol{1} = \boldsymbol{1}$, and $ \|\boldsymbol{A}-\frac{\boldsymbol{1}\boldsymbol{1}^\top}{N}\|< 1$.
\end{assumption}

Assumption \ref{Assum:matrix_A} ensures that the interaction graph induced by $\boldsymbol{A}$ is balanced and connected, i.e., there is a path from each agent to every other agent.

The optimization problem in [\ref{Equ:P1}] can be reformulated as the following multi-agent optimization problem:
\begin{equation}\label{Equ:P}
\begin{aligned}
\min_{\boldsymbol{x}\in \mathbb{R}^{N \times d}} f(\boldsymbol{x})=\frac{1}{N}\sum_{i=1}^{N}f_i(\boldsymbol{x}_i),\\
{\rm s.t.}\quad \boldsymbol{x}_1=\boldsymbol{x}_2=\cdots=\boldsymbol{x}_{N},
\end{aligned}
\end{equation}
where the matrix $\boldsymbol{x}$ is composed of all the local optimization variables, i.e., $\boldsymbol{x}=\left[ \boldsymbol{x}_1^\top;  \boldsymbol{x}_2^\top; \ldots; \boldsymbol{x}_N^\top \right] \in \mathbb{R}^{N \times d}$. \\
 In this paper, the local objective function $f_i(\boldsymbol{x}_i)$ and global objective function $f(\boldsymbol{x})$ can be nonconvex. They are assumed to satisfy the following conditions:
\begin{assumption}\label{Assum:Lipschitz}
Every $f_i(\cdot)$ is differentiable and is $\mathit{L_i}$-Lipschitz as well as $\rho_i$-Hessian Lipschitz: 
 \begin{equation}
    \begin{aligned}
        \| \nabla  f_i(\boldsymbol{x_1})-\nabla f_i(\boldsymbol{x_2}) \|
        \leqslant \mathit{L_i}\, \|\boldsymbol{x_1}-\boldsymbol{x_2} \|, \quad \forall \boldsymbol{x_1},\boldsymbol{x_2} \in \mathbb{R}^{d},
    \end{aligned}
 \end{equation}
\begin{equation}
     \begin{aligned}
            \|\nabla^2 f_i(\boldsymbol{x_1})-\nabla^2 f_i(\boldsymbol{x_2}) \|
            \leqslant \rho_i \, \|\boldsymbol{x_1}-\boldsymbol{x_2} \|, \quad \forall \boldsymbol{x_1},\boldsymbol{x_2} \in \mathbb{R}^{d}.
    \end{aligned}
\end{equation}

It can be verified that  the global gradient $\nabla F(\boldsymbol{\theta})=\frac{1}{N}\sum_{i=1}^{N}\nabla f_i(\boldsymbol{\theta})$ and Hessian $\nabla^2 F(\boldsymbol{\theta})=\frac{1}{N}\sum_{i=1}^{N}\nabla^2 f_i(\boldsymbol{\theta})$ are    $\mathit{L}$-Lipschitz and
$\rho$-Hessian Lipschitz, with $L=\frac{1}{N}\sum_{i} L_i$ and $\rho=\frac{1}{N}\sum_{i} \rho_i$.
\end{assumption}

As in most existing results on distributed nonconvex optimization, we assume that the local gradients  $\nabla f_i(\cdot)$ are bounded:

\begin{assumption}\label{Assum:Bound Gradient}
    There exists a constant G such that  $\|\nabla f_i(\boldsymbol{\theta})\| \leqslant G$ holds for all $\boldsymbol{\theta}\in\mathbb{R}^{d}$ and $i \in [N]$.
\end{assumption}

In this paper, we will show that quantization can help evade saddle points and ensure convergence to second-order stationary points in distributed nonconvex optimization. To this end, we first recall the following definitions for first-order stationary points, saddle points, and second-order stationary points, which are commonly used in the study of saddle-point problems:

\begin{definition}
    For a twice differentiable objective function $F(\cdot)$, we call $\boldsymbol{\theta^\star}\in\mathbb{R}^d$ a first-order (respt. second-order) stationary point if $\nabla F(\boldsymbol{\theta^\star})= \boldsymbol{0}$ (respt. $\nabla F(\boldsymbol{\theta^\star})= \boldsymbol{0}$ and $ \lambda_{\min}(\nabla^{2} F(\boldsymbol{\theta^\star} )) \geqslant 0 $) holds. Moreover, a first-order stationary point $\boldsymbol{\theta}^\star $ can be viewed as belonging to one of the three categories:
    \begin{itemize}
        \item   \emph{local minimum}: there exists a scalar $\gamma>0$ such that $F(\boldsymbol{\theta}^\star)\leqslant F(\boldsymbol{\theta})$ holds for any $\boldsymbol{\theta}$ satisfying $\|\boldsymbol{\theta}^\star-\boldsymbol{\theta}\|\leqslant\gamma$;
        \item   \emph{local maximum}: there exists a scalar $\gamma>0$ such that  $F(\boldsymbol{\theta}^\star)  \geqslant F(\boldsymbol{\theta})$ holds for any  $\boldsymbol{\theta}$ satisfying $\|\boldsymbol{\theta}^\star-\boldsymbol{\theta}\|\leqslant\gamma$;
        \item   \emph{saddle point}:  neither of the above two cases is true, i.e., for any scalar $\gamma>0$, there exist $\boldsymbol{\theta_1}$ and $\boldsymbol{\theta_2}$  satisfying $\|\boldsymbol{\theta_1}-\boldsymbol{\theta}^\star\|\leqslant \gamma$ and $\|\boldsymbol{\theta_2}-\boldsymbol{\theta}^\star\|\leqslant \gamma$ such that $F(\boldsymbol{\theta_1})<F(\boldsymbol{\theta}^\star)<F(\boldsymbol{\theta_2})$ holds.
    \end{itemize}
\end{definition}

Since distinguishing saddle points from local minima for smooth functions is NP-hard in general \cite{nesterov2000squared}, we focus on a subclass of saddle points, i.e., $\epsilon-$strict saddle points:

\begin{definition} {\rm ($\epsilon-$strict saddle point and $\epsilon-$second-order stationary point)}
\label{Def:strict saddle distributed}
For a twice-differentiable function $F(\cdot)$, we say that $\boldsymbol{\theta^\star} \in \mathbb{R}^d $ is an $\epsilon-$strict saddle point if 1) $\boldsymbol{\theta^\star}$ is an $\epsilon-$first-order
stationary point i.e., $\|\nabla F(\boldsymbol{\theta^\star})\| \leqslant \epsilon $; and 2) $\lambda_{\min}(\nabla^{2} F(\boldsymbol{\theta^\star} )) \leqslant -\sqrt{\rho\epsilon}$, where $\rho$ is the Hessian Lipschitz parameter in Assumption \ref{Assum:Lipschitz}. Similarly, $ \boldsymbol{\theta^\star} \in \mathbb{R}^d $ is an $\epsilon-$second-order stationary point if 1) $\boldsymbol{\theta^\star}$ is an $\epsilon-$first-order
stationary point, i.e., $\|\nabla F(\boldsymbol{\theta^\star})\| \leqslant \epsilon$ and 2) $\lambda_{\min}(\nabla^{2} F(\boldsymbol{\theta^\star} )) > -\sqrt{\rho\epsilon}$.
\end{definition}

 For a smooth function, a generic saddle point must satisfy that the minimum eigenvalue of its Hessian is non-positive. Our consideration of strict saddle points rules out the case where the minimum eigenvalue of the Hessian is zero.
 A line of recent work in the machine learning literature shows
that for many popular models in machine learning, all saddle points are indeed strict saddle points, with examples ranging from tensor decomposition \cite{ge2015escaping}, dictionary learning \cite{7755794}, smooth semidefinite programs \cite{boumal2016non}, to robust principal component analysis \cite{ge2017no}.

\section*{Proposed Algorithm}

By exploiting the effects of quantization, we propose a distributed nonconvex optimization algorithm that can ensure the avoidance of saddle points and convergence to a second-order stationary point. The detailed algorithm is summarized in Algorithm \ref{alg:1}.

\begin{algorithm}
    \caption{Distributed Optimization with Guaranteed Saddle-point Avoidance}
    \label{alg:1}
        \begin{itemize}[label={}]
            \item \textbf{Initialization:} $\boldsymbol{x}_{i}^{0} \in \mathbb{R}^d$ for every agent $i$;
            \item \textbf{Parameters:}
                            Stepsize sequences $\{\varepsilon_k\}$ and $\{\eta_k\}$;
            \item   \hspace{6em} Quantization level $\ell$;

            \item  \textbf{for}  {$k = 1,2,...$}  \textbf{do}
            \item  \hspace{1em} \textbf{for all}  {$i \in [N]$} \textbf{do}
            \item   \begin{itemize}[label={}]
                       \setlength{\leftmargin}{25.5em}

                       \item 1. Quantize its decision vector $\boldsymbol{x}_{i}^{k}$ to obtain $Q_{\ell}(\boldsymbol{x}_{i}^{k})$ and send the quantized $Q_{\ell}(\boldsymbol{x}_{i}^{k})$ to all neighbor agents in $\mathcal{N}_i$;

                       \item 2. Receive $Q_{\ell}(\boldsymbol{x}_{j}^{k})$  from neighbor agents $j \in \mathcal{N}_i$ and calculate the following estimate of the global optimization variable:
                            \begin{equation}
                            \begin{aligned}\label{iteration_i}
                                \tilde{\boldsymbol{x}}_i^{k+1}=\boldsymbol{x}_{i}^{k}+\varepsilon_k\sum_{j\in \mathcal{N}_i \cup \{ i \} }a_{ij}(Q_{\ell}(\boldsymbol{x}_{j}^{k})-\boldsymbol{x}_{i}^{k});
                                \end{aligned}
                            \end{equation}

                        \item 3. Calculate local gradient $\nabla
                            f_i(\boldsymbol{x}_{i}^{k})$  and update $\boldsymbol{x}_{i}^{k+1}$  by:
                                \begin{equation}
                                \begin{aligned}\label{iteration_i2}
                                 \boldsymbol{x}_{i}^{k+1}= \tilde{\boldsymbol{x}}_i^{k+1}-\eta_k\nabla f_i(\boldsymbol{x}_{i}^{k}).
                                \end{aligned}
                                \end{equation}
                   \end{itemize}
                \item  \hspace{1em} \textbf{end for}
                \item  \textbf{end for}
        \end{itemize}
\end{algorithm}



As key components of our approach to evading saddle points and ensuring convergence accuracy, we propose the following quantization scheme and stepsize strategy:

\subsection*{Quantization Scheme}

Our quantization scheme is inspired by the QSGD quantization scheme proposed in \cite{alistarh2017qsgd} and the TernGrad quantization scheme in \cite{wen2017terngrad}. (Note that the QSGD and TernGrad schemes were proposed to quantize gradients, whereas our Algorithm \ref{alg:1} quantizes optimization variables.) More specifically, at each time instant, we represent a continuous-valued variable with a randomized rounding to a set of quantization points with adjustable discrete quantization levels in a way that preserves the statistical properties of the original. However,  different from \cite{alistarh2017qsgd}, to ensure saddle-point avoidance, we employ two sets of quantization levels and purposely switch between the two sets of quantization levels in a periodic manner. The detailed scheme is described below:\\

 For any $ \boldsymbol{v} = [ v_1, v_2, \ldots, v_d ] \in \mathbb{R}^d$,
\begin{enumerate}
  \item   At any even-number iteration ($k$ is even), map every $v_i \in \mathbb{R}$ onto the quantization level-set: $\{\cdots,\, -3\ell,\, -2\ell,\, -\ell,\, 0,\, \ell,\, 2\ell,\, 3\ell,\,\cdots\}$ (which we will refer to ``level-set 1" hereafter) as follows:
  \begin{equation}\label{Q_s1}
    \begin{aligned}
   Q_\ell \left( v_i \right) =\begin{cases}
    	n\ell,& {\rm with\,\,probability}\,\,1-p\left( v_i, \ell \right)\\
    	(n+1)\ell,& {\rm with\,\,probability}\,\,p\left( v_i, \ell \right)\\
    \end{cases}
    \end{aligned}
    \end{equation}
    where $ n \in \mathbb{Z}$ is determined by the inequality $ n\ell \leqslant v_i <  (n+1)\ell$, and the probability $p   \left( v_i, \ell \right)$ is given by $p   \left( v_i, \ell \right)=\frac{v_i}{\ell}-n$.

  \item   At any odd-number iteration ($k$ is odd), map every $v_i \in \mathbb{R}$ onto the quantization level-set: $\{\cdots,\, -2.5 \ell,\, -1.5\ell,\, -0.5 \ell,\, 0.5 \ell,\, 1.5 \ell,\, 2.5\ell,\, \cdots  \}$ (which we will refer to ``level-set 2" hereafter) as follows:
  \begin{equation}\label{Q_s2}
    \begin{aligned}
   Q_\ell  \left( v_i \right) =\begin{cases}
    	(n'-0.5)\ell,& {\rm with\,\,probability}\,\,1-p'\left( v_i, \ell \right)\\
    	(n'+0.5)\ell,& {\rm with\,\,probability}\,\,p'\left( v_i, \ell \right)\\
    \end{cases}
    \end{aligned}
    \end{equation}
   where $ n' \in \mathbb{Z}$ is determined by the inequality $ (n'-0.5)\ell \leqslant v_i <  (n'+0.5)\ell$, and the probability $p'\left( v_i, \ell \right)$ is given by $p'\left( v_i, \ell \right)=\frac{v_i}{\ell}-n'+0.5$.
\end{enumerate}

It is worth noting that compared with existing quantization schemes such as \cite{alistarh2017qsgd}, this periodic switching between two sets of quantization levels does not introduce extra communication overheads. However, it avoids the possibility that any quantization input $v_i$ always coincides with an endpoint of a quantization interval, resulting in a deterministic quantization output. Namely, for any point $ \boldsymbol{v} \in \mathbb{R}^d$, it gives two different representations in the quantized space, which is key to perturb and avoid the state from staying on undesired saddle points under a non-zero stepsize (which will be elaborated later).

An instantiation of this quantization scheme is depicted in Fig.~\ref{Fig:Quantization}. It can be verified that the proposed quantization scheme satisfies the following properties:

\begin{figure}
  \centering
  \includegraphics[width=0.8\linewidth]{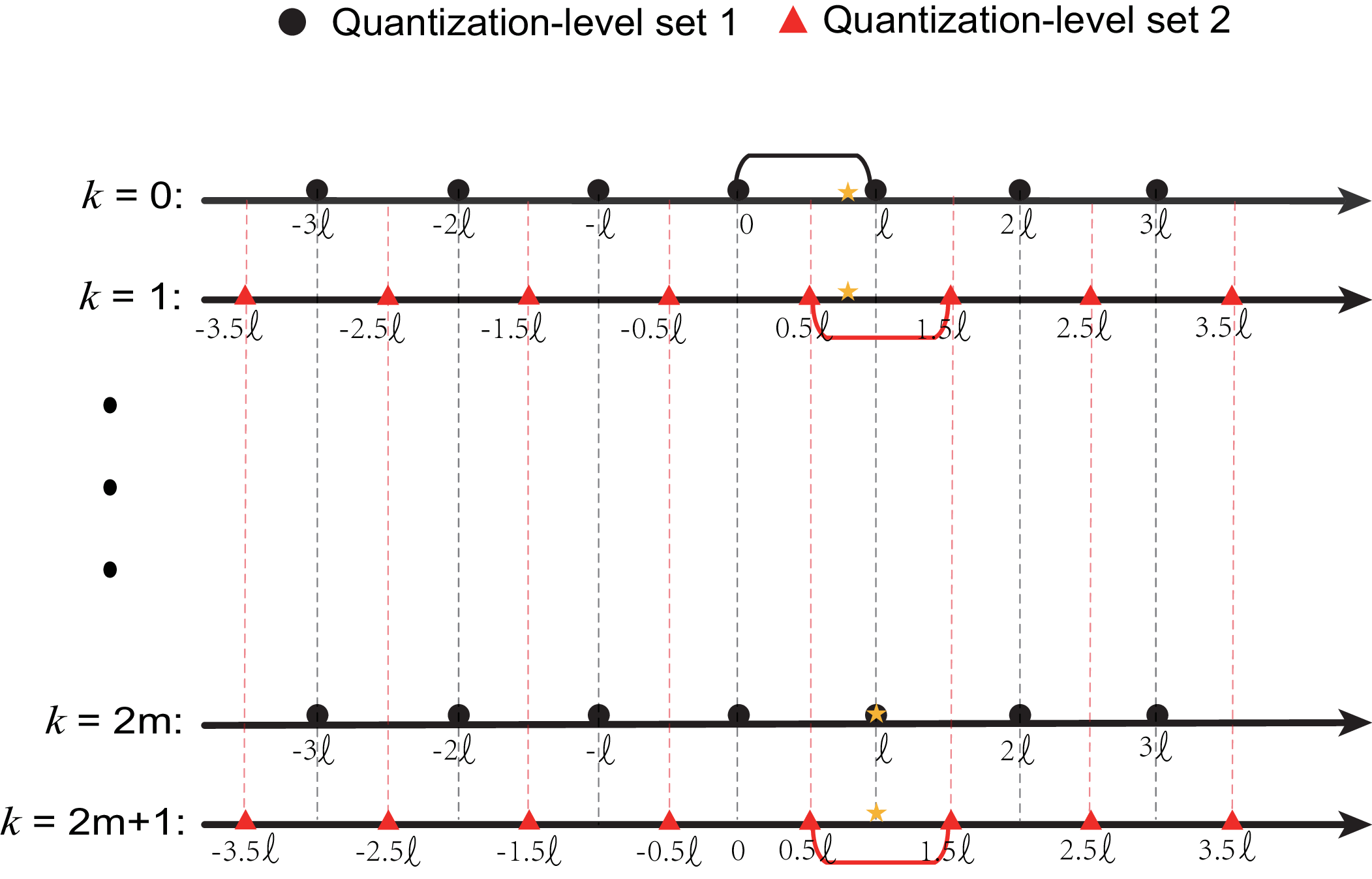}\\
  \caption{ The proposed quantization scheme with quantization interval $ \ell $. The star represents a value to be quantized, and it is located in the quantization interval of $[0, \ell]$ under level-set 1 and $[0.5\ell, 1.5 \ell]$ under level-set 2.
  At any even-number iteration ($k$ is even), the star value will be quantized to either 0 or $\ell$, with respective probabilities provided in [\ref{Q_s1}]. At any odd-number iteration ($k$ is odd), the star value will be quantized to either $0.5\ell$ or $1.5 \ell$, with respective probabilities given in [\ref{Q_s2}].}
  \label{Fig:Quantization}
\end{figure}



\begin{lemma}\label{lem:quantization}
For any   $  \boldsymbol{v}  \in \mathbb{R}^d $, our quantization scheme $ Q_\ell(\boldsymbol{v})=[ Q_\ell  \left( v_1 \right), Q_\ell  \left( v_2 \right),...,Q_\ell  \left( v_d \right)]$ has the following properties:
\begin{enumerate}
  \item Unbiased quantization: $\mathbb{E} \left[ {Q_\ell( \boldsymbol{v} )}  \right] =  \boldsymbol{v} $,
  \item Bounded variance: $\mathbb{E} \left[ { \|Q_\ell ( \boldsymbol{v} )- \boldsymbol{v} \|^2 }  \right] \leqslant d \ell^2 $.
\end{enumerate}
\end{lemma}


\subsection*{Stepsize Strategy}

 In addition to purposely employing switching in the quantization scheme, the stepsizes $\{\varepsilon_k\}$ and $\{\eta_k\}$ in Algorithm \ref{alg:1} also have to be judiciously designed so as to evade saddle points and ensure convergence to a second-order stationary point. Intuitively speaking, in the early stage where saddle points may trap the optimization process, the stepsize $ \varepsilon_k $ should be large enough to ensure that the switching quantization-induced perturbation can effectively stir the evolution of optimization variables. However, to ensure that the optimization process can converge to an optimal solution, the quantization effect should gradually diminish, or in other words, $ \varepsilon_k $ should converge to zero. In addition, in distributed optimization, to ensure that all agents can converge to an optimal solution without any error, the stepsize $ \eta_k $ also has to converge to zero (different from the centralized case, in distributed optimization, a constant stepsize will lead to optimization errors that are in the order of the stepsize \cite{nedic2009distributed, vlaski2021distributed2, reisizadeh2019exact}). Moreover, to ensure that all agents can converge to the same optimal solution, the stepsize $ \varepsilon_k $ should decay slower than $ \eta_k $ \cite{wang2023tailoring, wang2022ensure, wang2022decentralized, doan2020convergence}. To fulfill these requirements, we design the stepsize sequences $\{\varepsilon_k\}$ and $\{\eta_
k\}$ as follows:

\begin{enumerate}
    \item Choose two positive constants $\alpha$ and $\beta$ sequentially that satisfy the following relations: $0.6< \alpha <\frac{2}{3}$ and $\frac{3}{2} \alpha  < \beta < 1$.
    And then use these constants to construct two reference functions $ \frac{c_1}{1+c_2 t^\alpha }$ and $\frac{c_1}{1+c_2 t^\beta }$,
    where $t$ is continuous time and $c_1$ and $c_2$ are all positive constants.
\\
    \item For any probability $ p $ (where $1-p$ represents the desired probability of converging to a second-order stationary point, which can be chosen to be arbitrarily close to one, see the statement of Theorem \ref{Th:Second-order} for details) and $\epsilon >0$ given in Definition \ref{Def:strict saddle distributed}, select: \\
    $t_0 \geqslant \max \{ {C}_1, {C}_2, {C}_3 \}$, $t_{i+1}=t_{i}+ \lceil \frac{1+c_2 t_i^\alpha }{c_1\sqrt{\rho\epsilon}} \rceil $
    for $1\leqslant i\leqslant I$,
    where  $I= 30 \max \{ \frac{f_0 -f^\star}{Q},\frac{2(f_0 -f^\star)\varepsilon _{t_0}}{\epsilon ^2 \eta _{t_0}}\} $.\footnote[1]{ $  {C}_1=( \frac{4c_{1}^{2/3}(d_1+d_2)}{pc_{2}^{2/3}( 1-\sigma _2 )} ) ^{\frac{3}{2\alpha}}$, $  {C}_2=( \frac{4( f_0-f^{\star} ) ( d_1+d_2 ) ^{2/3}( 1-\sigma _2 ) ^{2/3}c_1}{c_2p\epsilon ^2\sqrt{\rho \epsilon}} ) ^{\frac{1}{2\alpha -\beta}}$, ${C}_3=(\frac{12\rho (d_1+d_2)^{1/6}}{(1-\sigma _2)^{1/6}\sqrt{\gamma}(\rho \epsilon )^{1/4}\ell})^{\frac{1}{\beta -4\alpha /3}}$,  $Q = \frac{1}{60^2}\sqrt{\frac{\epsilon ^3}{\rho}}$, where $d_1=\frac{1+( 1-\sigma _2 ) \varepsilon _0}{1-\sigma _2}G^2$, $d_2=( 1+( 1-\sigma _2 ) \varepsilon _0 ) {\sigma _2}^2Nd\ell^2$, and $\sigma _2$ is the second largest eigenvalue of $\boldsymbol{A}$.
    $f_0$ is the objective function value at $k=0$.
    $f^\star$ denotes an estimated lower bound on the minimum global objective function $f(\cdot)$. For instance, in the matrix factorization problem where a low-rank matrix $U \in \mathbb{R} ^{d \times r}$ is used to approximate a high-dimension matrix $ \boldsymbol{M^\star} \in \mathbb{R} ^{d \times d}$, the objective function is  $ f(\boldsymbol{U}) = \frac{1}{2} \|\boldsymbol{U}\boldsymbol{U}^\top -\boldsymbol{M^\star} \|_{F}^{2} $ and we can use $f^\star = 0$ as the lower bound \cite{jin2017escape}.
    }
\\

    \item  The sequences $\{\varepsilon_k\}$ and $\{\eta_k\}$  for $\forall k \in \mathbb{Z^{+}} $  are given as follows:\\
    \begin{equation}\label{Equ:varepsilon}
        \begin{aligned}
           \varepsilon _k=\begin{cases}
	\frac{c_1}{1+c_2k^{\alpha}},&		k <  t_0\\
	\frac{c_1}{1+c_2{t_i^{\alpha}}},&		t_i\leqslant k  < t_{i+1}\\
	\frac{c_1}{1+c_2k^{\alpha}},&		k\geqslant t_I\\
        \end{cases}
        \end{aligned}
    \end{equation}
    \begin{equation}\label{Equ:eta}
        \begin{aligned}
            \eta _k=\begin{cases}
        	\frac{c_1}{1+c_2k^{\beta}},&		k < t_0\\
        	\frac{c_1}{1+c_2{t_i^{\beta}}},&		t_i\leqslant k  < t_{i+1}\\
        	\frac{c_1}{1+c_2k^{\beta}},&		k\geqslant t_I\\
                     \end{cases}
        \end{aligned}
    \end{equation}
\end{enumerate}

\begin{figure}
    \centering
    \includegraphics[width=1.1\linewidth]{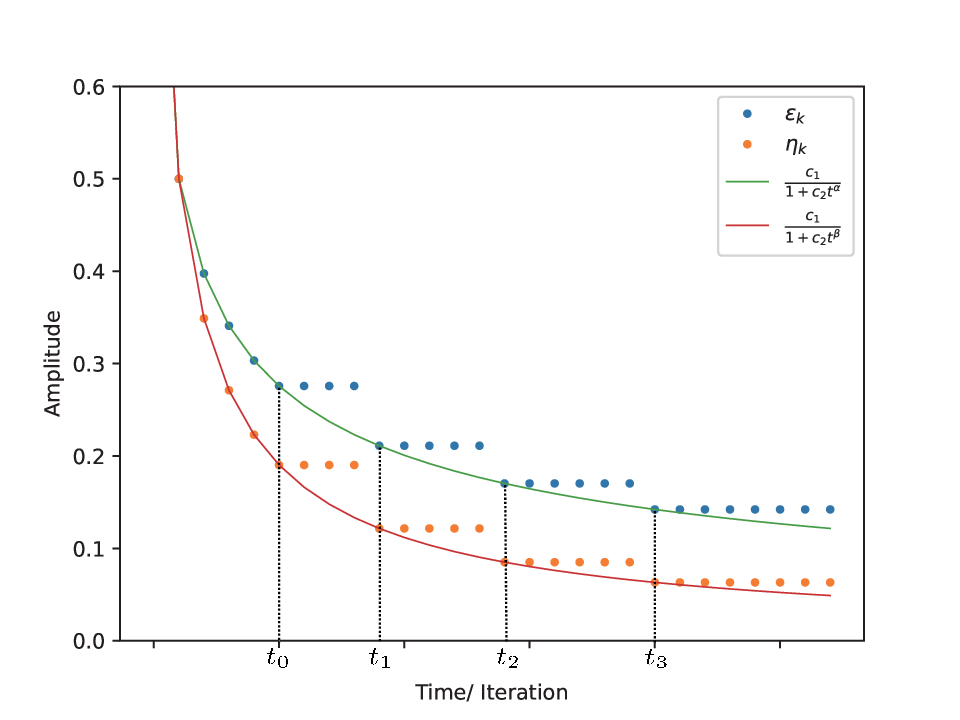}
    \caption{An illustrative example of the stepsizes. The two solid curves represent two reference functions which are defined on the continuous time $t$. The blue and orange dots represent the values of stepsizes $\varepsilon_k$ and $\eta_k$ at discrete time instants $k$ (which are periodic samples of the continuous time $t$).  The time instants $t_0, t_1, t_2, t_3 $ are determined in the second step of the stepsize strategy. Before $t_0$, the descent of the stepsize sequences is aligned with the reference functions. In intervals $[t_i, t_{i+1})$, the stepsizes remain constant, as described in the third step of the stepsize strategy.}
    \label{fig:stepsize}
\end{figure}

The setup in [\ref{Equ:varepsilon}] and [\ref{Equ:eta}] makes the stepsizes $ \varepsilon_k $ and $ \eta_k $ follow a decrease-and-hold pattern, as illustrated in Fig.~2.
The rationale for this design can be understood intuitively as follows: To ensure that all agents can converge to a desired optimal solution, the iteration process must fulfill two objectives simultaneously: 1) ensure the consensual convergence of all agents to a stationary point; and 2) avoid saddle points. The ``decrease" stages are important to fulfill the first objective in that    the
stepsizes  $\varepsilon_k$ and $\eta_k$ are used to attenuate the quantization error and input heterogeneity among the agents, respectively, (both of which act as counter-forces for reaching consensus among all agents' iterates), and hence their decrease is key to ensure reaching consensus among the agents.
 The ``hold" stages are necessary to accumulate enough stochastic quantization effects to stir the
evolution of optimization variables and ensure saddle-point avoidance. It is worth noting that the decrease of stepsize $\eta_k$ should also be carefully designed in both decreasing speed and timespan to ensure that sufficient gradient descent can be carried out to explore the solution space and ensure convergence to a stationary point. Hence we judiciously design the stepsize strategy to strike a balance between accumulating quantization noise to evade saddle points and attenuating quantization noise to ensure consensual convergence of all agents to a desired stationary point.

Based on the proposed quantization scheme and stepsize strategy, we can prove that the proposed Algorithm \ref{alg:1} can ensure all agents to evade saddle points and converge to the same second-order stationary point. For convenience of exposition, we divide the convergence analysis into two parts: ``Consensual convergence to a first-order stationary point" and ``Escaping saddle points and converging to a second-order stationary point". We leave all proofs in the Supporting Information.

\section*{Consensual Convergence to a First-order Stationary Point}\label{Sec:consensus}
We first prove that the proposed algorithm can ensure all agents to reach consensus on their optimization variables. For the convenience of analysis, we represent the effect of quantifying $ \boldsymbol{x}_i^k$ as adding noise to $ \boldsymbol{x}_i^k$, i.e., $Q_{\ell}(\boldsymbol{x}_i^k)=\boldsymbol{x}_i^k+\boldsymbol{\xi}_i^k$, where $\boldsymbol{\xi}_i^k$ is the stochastic quantization error.
Using the iteration dynamics in [\ref{iteration_i}] and [\ref{iteration_i2}], we can obtain the following relationship:
\begin{equation}\label{Equ:x_i-iteration}
\begin{aligned}
 \boldsymbol{x}^{k+1}_i
&=(1-\varepsilon _k) \boldsymbol{x}_i^k+\varepsilon _k\boldsymbol{A} Q_{\ell}(\boldsymbol{x}_i^k)-\eta _k \nabla f_i(\boldsymbol{x}_i^k).
\end{aligned}
\end{equation}

By defining $ \boldsymbol{x}^k=\left[ (\boldsymbol{x}_1^k)^\top;  (\boldsymbol{x}_2^k)^\top; \cdots; (\boldsymbol{x}_N^k)^\top \right] \in \mathbb{R}^{N \times d}$, $\boldsymbol{A}_k=(1-\varepsilon_k)\boldsymbol{I}+\varepsilon_k \boldsymbol{A}$, $\nabla f( \boldsymbol{x}^{k} )=\left[ \nabla f_1^\top(\boldsymbol{x}_1^k) ; \nabla f_2^\top(\boldsymbol{x}_2^k); \cdots ; \nabla f_N^\top(\boldsymbol{x}_N^k) \right] \in \mathbb{R}^{N \times d}$, and $\boldsymbol{\xi}^k=\left[
	(\boldsymbol{ \boldsymbol{\xi}_1^k})^\top ; (\boldsymbol{\xi}_2^k)^\top; \cdots ; (\boldsymbol{\xi}_N^k)^\top \right] \in \mathbb{R}^{N \times d}$, we can recast the relationship in [\ref{Equ:x_i-iteration}] into the following more compact form:

\begin{equation}\label{iteration}
\begin{aligned}
      \boldsymbol{x }^{k+1}  &= \boldsymbol{A}_k\boldsymbol{x}^{k} +\varepsilon_k \boldsymbol{A} \boldsymbol{\xi}^k-\eta_k\nabla f(\boldsymbol{x}^{k}).
\end{aligned}
\end{equation}

Let $\boldsymbol{\bar{x}}^k$ be the average of all local optimization variables, i.e., $ \boldsymbol{\bar{x}}^k = \frac{1}{N}\sum_{i=1}^N \boldsymbol{x}_i^k$. It can be verified that $\boldsymbol{\bar{x}}^k$ is equal to $ \frac{({  \boldsymbol{x}^{k}  })^\top \boldsymbol{1}}{N}  $, which can be further verified to satisfy the following relationship based on [\ref{iteration}]:
\begin{equation}
    \begin{aligned}
     \label{y^k}
\boldsymbol{\bar{x}}^{k+1}=\boldsymbol{\bar{x}}^k+\varepsilon_k \frac{ ( \boldsymbol{\xi}^k  )^\top \boldsymbol{1}}{N}- \eta_k\frac{\nabla f^\top(\boldsymbol{x}^{k} )\boldsymbol{1}}{N}.
    \end{aligned}
\end{equation}

Define the consensus error between individual agents' local optimization variables and the average optimization variable $\boldsymbol{\bar{x}}^k$ as $\boldsymbol{e}^k:=\boldsymbol{x}^{k}-\boldsymbol{1}(\boldsymbol{{\bar{x}}}^{k})^\top$. It can be verified that the $i$-th row of $\boldsymbol{e}^k$, i.e., $\boldsymbol{e}_i^k$, satisfies $\boldsymbol{e}_i^k=(\boldsymbol{x}_i^k)^\top -(\boldsymbol{{\bar{x}}}^{k})^\top$.   Using the algorithm iteration rule described in [\ref{iteration_i}] and [\ref{iteration_i2}], we can obtain the following iteration dynamics for $\boldsymbol{e}^k$:
\begin{equation}\label{Equ:e^k1}
     \begin{aligned}
    \boldsymbol{e}^{k+1}
           =&\boldsymbol{A}_k\boldsymbol{e}^k+\varepsilon _k\boldsymbol{A}\boldsymbol{W}\boldsymbol{\xi}^k-\eta _k\boldsymbol{W}\nabla f(\boldsymbol{x}^{k}),
    \end{aligned}
\end{equation}
where $\boldsymbol{W}=\boldsymbol{I}-\frac{\boldsymbol{1}\boldsymbol{1}^{\top}}{N}$.

Based on the dynamics of consensus errors $\boldsymbol{e}^k$ in [\ref{Equ:e^k1}], we can prove that the consensus error $\|\boldsymbol{e}^k\|^2$ will converge almost surely to zero, i.e., all  $\boldsymbol{x}_i^k$ will almost surely converge to the same value.


\begin{theorem} {\rm (Consensus of Optimization Variables)} \label{corollaryk_0}
  Let Assumptions \ref{Assum:matrix_A}, \ref{Assum:Lipschitz}, and \ref{Assum:Bound Gradient} hold. Given any probability $0<p<1$, Algorithm \ref{alg:1} with our stepsize strategy (which takes $p$ as input) ensures consensus error less than $ \mathcal{O} \left(\frac{1}{k}\right)^{\frac{\alpha}{3}}$  with probability at least $1-p$ for all $k \geqslant t_0$, where $t_0$ is given in step 1 and step 2 of the stepsize strategy, respectively:
    \begin{equation}\label{Equ:consensus_rate_alpha}
    \begin{aligned}
       \mathbb{P} \left(\left\| \boldsymbol{e}^k \right\| ^2\leqslant \mathcal{O} \left(\frac{1}{k}\right)^{\frac{\alpha}{3}}, \,\,for\,\,all\,\,k\geqslant t_0\right)\geqslant 1-p.
    \end{aligned}
    \end{equation}
Moreover, all agents' optimization variables converge to the same value almost surely, i.e., the consensus error $\|\boldsymbol{e}^k\|$ converges almost surely to zero.

\end{theorem}



Based on the consensus result in Theorem \ref{corollaryk_0},  we can further prove that Algorithm \ref{alg:1} ensures all local optimization variables to converge to a first-order stationary point under the given quantization scheme and stepsize strategy:

\begin{theorem}{\rm (Converging to a First-order Stationary Point)}\label{Th:first-order}
Let Assumptions \ref{Assum:matrix_A}, \ref{Assum:Lipschitz}, and \ref{Assum:Bound Gradient} hold. Given any probability $0<p<1$, Algorithm \ref{alg:1} with our stepsize strategy
 (which takes $p$ as input) ensures that the gradient $\|  \nabla F\left( \boldsymbol{\bar{x}}^k \right) \|$ will converge to zero with a probability no less than $1-p$, i.e.,
\begin{equation}\label{Equ:first-order}
  \mathbb{P} \left( \lim_{k\rightarrow \infty}\left\| \nabla F\left( \boldsymbol{\bar{x}}^k \right) \right\| ^2 = 0 \right) \geqslant 1-p.
\end{equation}
\end{theorem}
It is worth noting that due to the employment of $\varepsilon_k$ (which gradually suppresses the influence of quantization errors) and the unbiasedness of the quantization scheme (the quantization error has a mathematical
expectation equal to zero),  our algorithm ensures convergence to an exact minimum that has a zero gradient value (with zero steady-state error). In fact, the absence of steady-state error under unbiased quantization has been obtained in the literature such as QSGD \cite{alistarh2017qsgd} and TernGrad \cite{wen2017terngrad}.

\section*{Escaping Saddle Points and Converging to a Second-order Stationary Point}




According to Definition \ref{Def:strict saddle distributed}, saddle points are undesirable states that stall the iteration process. Given that 1) individual agents' local optimization variables $\boldsymbol{x}_i^k$ quickly converge to the same value (reach consensus) according to Theorem \ref{corollaryk_0}; and 2) before reaching consensus, inter-agent interaction acts as an additional force (besides the gradient) to keep individual states $\boldsymbol{x}_i^k$  evolving and hence to avoid them from being trapped at any fixed value, we can only consider the saddle-point problem when the states are consensual. In fact, even after all states have reached consensus, since the force brought by inter-agent iteration diminishes at a slower rate than the driven force of the gradient ($\varepsilon_k$ decays slower than $ \eta_k $ in our stepsize strategy), the quantized interaction will have enough perturbations on individual agents' optimization variables to efficiently avoid them from being trapped at any saddle point. Formally, we can prove the following results:

\begin{theorem}[Escaping Saddle Points]\label{Th:Escaping Saddle Points}

Let Assumptions \ref{Assum:matrix_A}, \ref{Assum:Lipschitz}, and \ref{Assum:Bound Gradient} hold. Given any probability $0<p<1$, Algorithm \ref{alg:1} with our stepsize strategy
 (which takes $p$ as input) ensures that any ``holding stage" in the stepsize strategy reduces the objective function by a substantial amount. More specifically,  for any $i \in \{1,2,\ldots I\}$,  after no more than $K=\mathcal{O}(\frac{1}{\varepsilon_{t_i}\sqrt{\rho\epsilon}})$ iterations with the stepsizes held at  $\{\varepsilon_{t_i}, \eta_{t_i}\}$, Algorithm \ref{alg:1} ensures that with a substantial probability, the objective function has a significant decrease, i.e.,
\begin{equation}
\mathbb{P} \left( F\left( \boldsymbol{\bar{x}}^{t_i+K} \right) -F\left( \boldsymbol{\bar{x}}^{t_i} \right) \leqslant -Q \right) \geqslant \frac{1}{3}-p,
\end{equation}
where $Q$ is a constant satisfying $Q=\mathcal{O}\left(  \sqrt{\frac{\epsilon ^3}{\rho}} \right)$.

\end{theorem}


It is worth noting that although the inter-agent interaction (after quantization) can perturb individual agents' optimization variables from staying at any fixed point in the state space, it cannot ensure escaping from a saddle point since the state may evolve in and out of the neighborhood of a saddle point. To facilitate escaping from saddle points, we have to make full use of the existence of descending directions at strict saddle points. More specifically, in our design of the quantization scheme and stepsize strategy, we exploit random quantization to ensure that perturbations exist in every direction and use switching quantization levels to ensure that the amplitude of such perturbations is persistent. To ensure a sufficient integration of the perturbation effect into the iterative dynamics and make it last long enough to evade a saddle point, we hold the stepsizes $\varepsilon_k$ and $\eta_k$ constant for a judiciously calculated period of time (see Fig.~\ref{fig:stepsize}).

In fact, besides evading a saddle point, Theorem \ref{Th:Escaping Saddle Points} establishes that in each ``holding stage" where the stepsizes $\varepsilon_k$ and $\eta_k$ are held constant, the algorithm is guaranteed to decrease in the function value for a significant amount. Therefore, if we can have an estimation of a lower bound on the optimal function value $f^\star$, we can repeat this holding stage multiple times to ensure avoidance of all potentially encountered saddle points, and hence, to ensure convergence to a second-order stationary point.

 In practice, during the algorithm's iterations, encountered points can be classified into two categories: points with relatively large gradients $\|\nabla F(
 \boldsymbol{\bar{x}} )\| > \epsilon $ and points with small gradients $\|\nabla F( \boldsymbol{\bar{x}} )\| \leqslant \epsilon $, i.e., saddle points. We can prove that within the $t_I$ iterations defined in the stepsize strategy, the algorithm will encounter a second-order stationary point at least once:
\begin{theorem}{\rm (Converging to a Second-order Stationary Point)}\label{Th:Second-order}
Let Assumptions \ref{Assum:matrix_A}, \ref{Assum:Lipschitz}, and \ref{Assum:Bound Gradient} hold. For any $\epsilon>0$ and any given probability $0<p<1$, our stepsize strategy (which takes $p$ as input) ensures that Algorithm \ref{alg:1} will visit an $\epsilon-$second-order stationary point at least once with probability at least $1-p$ in $t_I$ iterations stated in the stepsize strategy.
\end{theorem}

From the derivation of Theorem \ref{Th:Second-order} in the Supporting Information, we can obtain that it takes the following number of  iterations to find an $\epsilon-$second-order stationary point:
\begin{equation}
\begin{aligned}
\mathcal{O} \left(  \frac{1}{\epsilon ^2}\max \left\{
  (  Nd\ell ^2 ) ^{\frac{3}{2\alpha}}, ( \frac{( Nd\ell ^2 ) ^{\frac{2}{3}}}{\epsilon ^{2.5}} )^{\frac{1}{2\alpha -\beta}}, (\frac{\left( Nd \right) ^{1/6}}{\epsilon ^{1/4}\ell ^{2/3}} ) ^{\frac{1}{\beta -4\alpha /3}}
\right\} \right),
\end{aligned}
\end{equation}
where  $N$ is the  number of agents participating in the distributed optimization, $d$ is the dimension of the optimization variable,    $\ell$ is the size of the quantization interval, and $\alpha$ and $\beta$ are the parameters in stepsizes $\varepsilon_k$ and $\eta_k$, respectively (note that $2\alpha>\beta>\frac{4}{3}\alpha$ holds according to our stepsize strategy). Therefore, the computational complexity of our algorithm increases polynomially with increases in the network size $N$ and the dimension of optimization variation  $d$.  It is worth noting that the computational complexity does not increase monotonically with the size of the quantization interval $\ell$: both a too small $\ell$ and a too large $\ell$ lead to a high computational complexity.   This is understandable since a too small quantization interval $\ell$ leads to too small quantization errors to stir the
evolution of optimization variables, which makes it hard to evade saddle points; whereas a too large quantization interval $\ell$ results in too much noise injected into the system, which is also detrimental to the convergence of all agents to a stationary point.

\section*{Experiments}
\begin{figure}
\centering
\includegraphics[width=.6\linewidth]{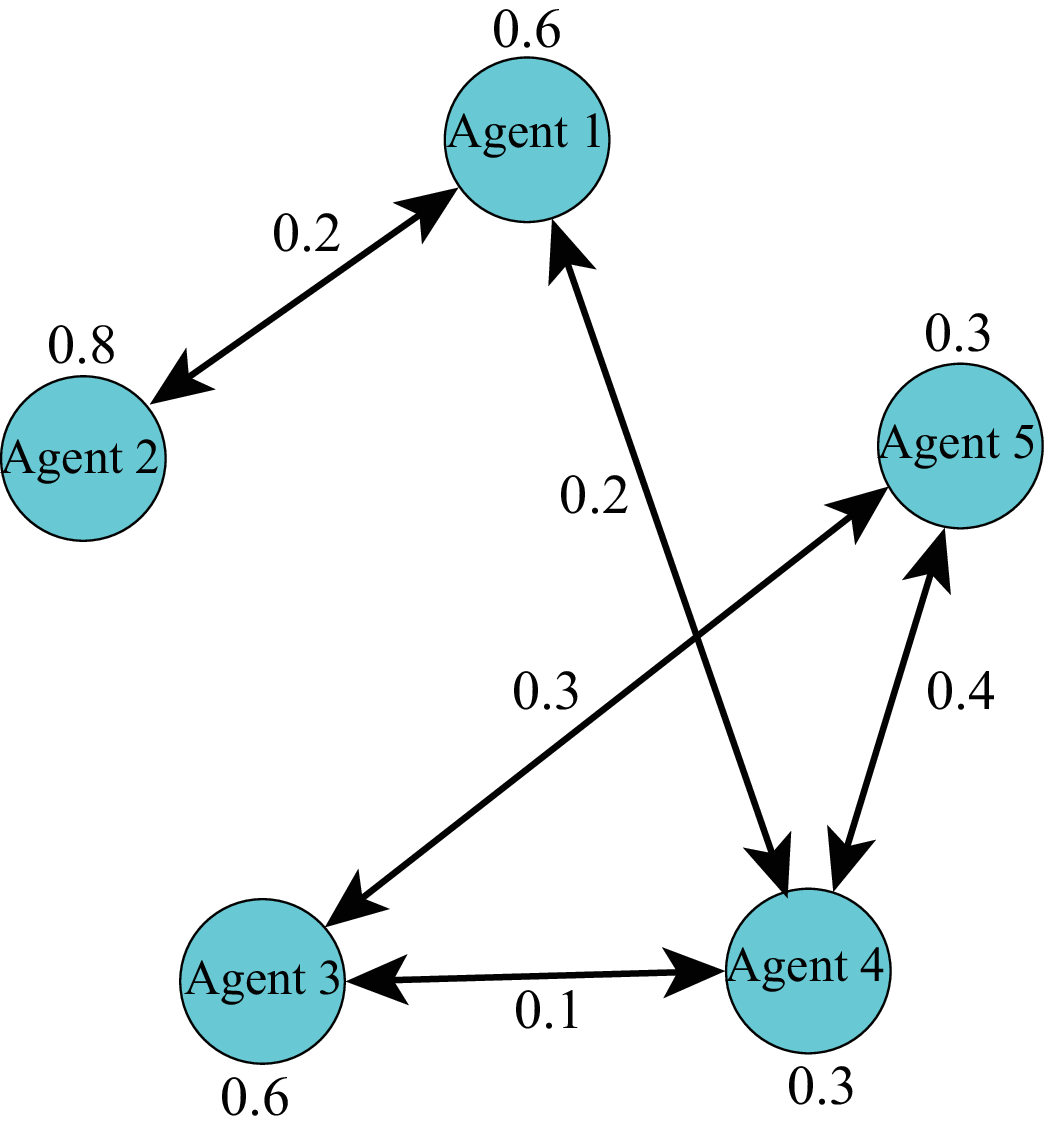}
\caption{Interaction weights of five agents}
\label{Fig:weights}
\end{figure}

In this section, we evaluate the performance of the proposed algorithm in five nonconvex-optimization application examples with different scales and complexities. In all five experiments, we consider five agents interacting on the topology depicted in Fig.~\ref{Fig:weights}. \\

\subsection*{Binary Classification}
In this experiment, we consider a simple $\{0, 1\}$ -- classification neural network with a single linear hidden layer and a logistic activation function. We use the cross-entropy loss function to train the network (see \cite{vlaski2021distributed} for details).
We denote the feature vector as $\boldsymbol{h} \in \mathbb{R}^M$ and the binary class label as $y\in \{-1, 1\}$. For the fully connected hidden layer, we represent the weights as  $\boldsymbol{W}_2 \in \mathbb{R}^{ L \times M }$ and $\boldsymbol{W}_1 \in \mathbb{R}^L $. The output is of the form:
\begin{equation}
\hat{y} = \frac{1}{1+e^{-\left< \boldsymbol{h}, \boldsymbol{W}_2^{\top}\boldsymbol{W}_1 \right>}}
\end{equation}
Under the commonly used cross-entry loss function, the objective function is of the following form:
\begin{equation}
L\left( \boldsymbol{W}_1,\boldsymbol{W}_2 \right) =\log \left( 1+e^{-y\left< \boldsymbol{h},\boldsymbol{W}_2^{\top}\boldsymbol{W}_1 \right>} \right)
\end{equation}

\begin{figure}
\centering
\includegraphics[width=.8\linewidth]{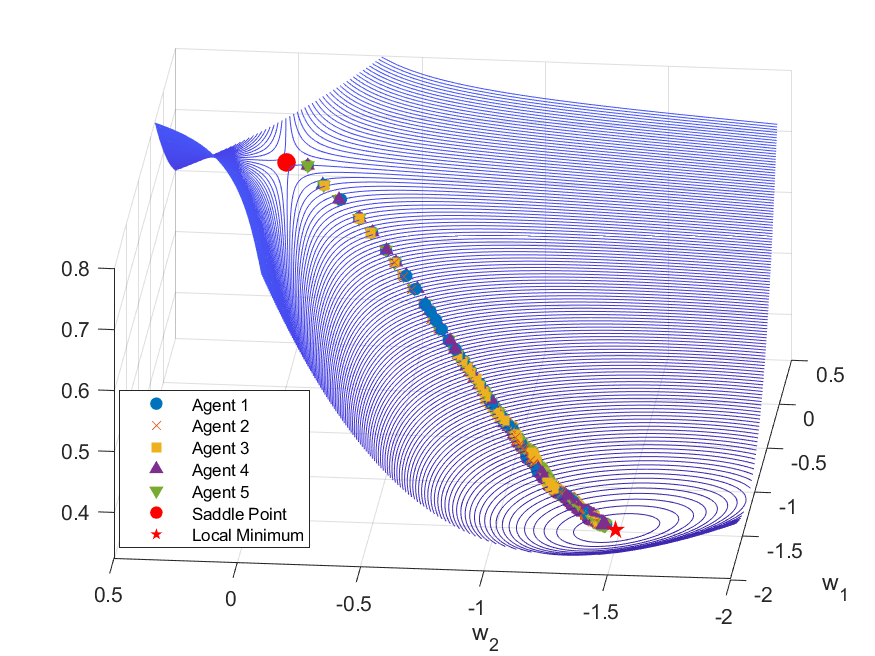}
\caption{Trajectories of all five agents when initialized on the saddle point (0,0). Note that all trajectories overlap with each other, implying perfect consensus among the agents.}
\label{Fig: trajectories}
\end{figure}

To visualize the evolution of optimization variables under our algorithm, we consider the scalar case with $L=M=1$ and plot the expected loss function (with regulation) in Fig.~\ref{Fig: trajectories}:
\begin{equation}\label{Equ:Logistic}
  F(w_1, w_2)=\mathbb{E} \left[ L\left( w_1, w_2 \right) \right]+\frac{\rho}{2}\left( \parallel w_1\parallel ^2+\parallel w_2\parallel ^2 \right)
\end{equation}

When the training samples satisfy $\mathbb{E} \left[ y\boldsymbol{h}\right]=1$ and the regularization parameter is set to $\rho = 0.1$, it becomes apparent that $(w_1, w_2)=(0,0)$ is a saddle point. We can also verify that this saddle point is a strict saddle point since its Hessian has a negative eigenvalue of $-0.4$.
In our numerical experiment, we purposely initialize all the agents from the strict saddle point $(0, 0)$, and plot in Fig.~\ref{Fig: trajectories} the evolution of each agent under stepsize parameters $\alpha= 0.62$, $\beta = 0.94$, $c_1 = 0.03$, and $c_2 = 0.3$. It can be seen that due to the quantization effect, all five agents collectively move along the descending direction, implying that our algorithm can effectively evade saddle points.

\subsection*{Matrix Factorization}
In this experiment, we consider the `Matrix Factorization' problem using the `MovieLen 100K' dataset and compare the performance of the proposed algorithm with a commonly used algorithm in \cite{yuan2016convergence}. In the matrix factorization problem, given a matrix $\boldsymbol{A}\in \mathbb{R} ^{m\times n}$ and $r < \min \{m, n\}$, the goal is to find two matrices $\boldsymbol{U} \in \mathbb{R} ^{m\times r}$ and $\boldsymbol{V} \in \mathbb{R} ^{n\times r}$ such that $F(\boldsymbol{U}, \boldsymbol{V}) = \frac{\|\boldsymbol{U}\boldsymbol{V}^\top -\boldsymbol{A}\|_{F}^{2}}{2}$ is
minimized.
However, due to the invariance property \cite{zhu2021global}, the matrix factorization problem cannot be considered strongly convex (or even convex) in any local neighborhood around its minima. In our numerical experiments, we implement both our algorithm and the algorithm in \cite{yuan2016convergence}.  In order to ensure a fair comparison, both algorithms share the same set of learning rates ($\alpha= 0.62$, $\beta = 0.94$, $c_1 = 0.3$, $c_2 = 0.3$). For the quantization scheme, we chose $\ell$ such that all quantized outputs are representable using a binary string of 9 bits. We spread the data evenly across the five agents.

\begin{figure}
\centering
\includegraphics[width=.8\linewidth]{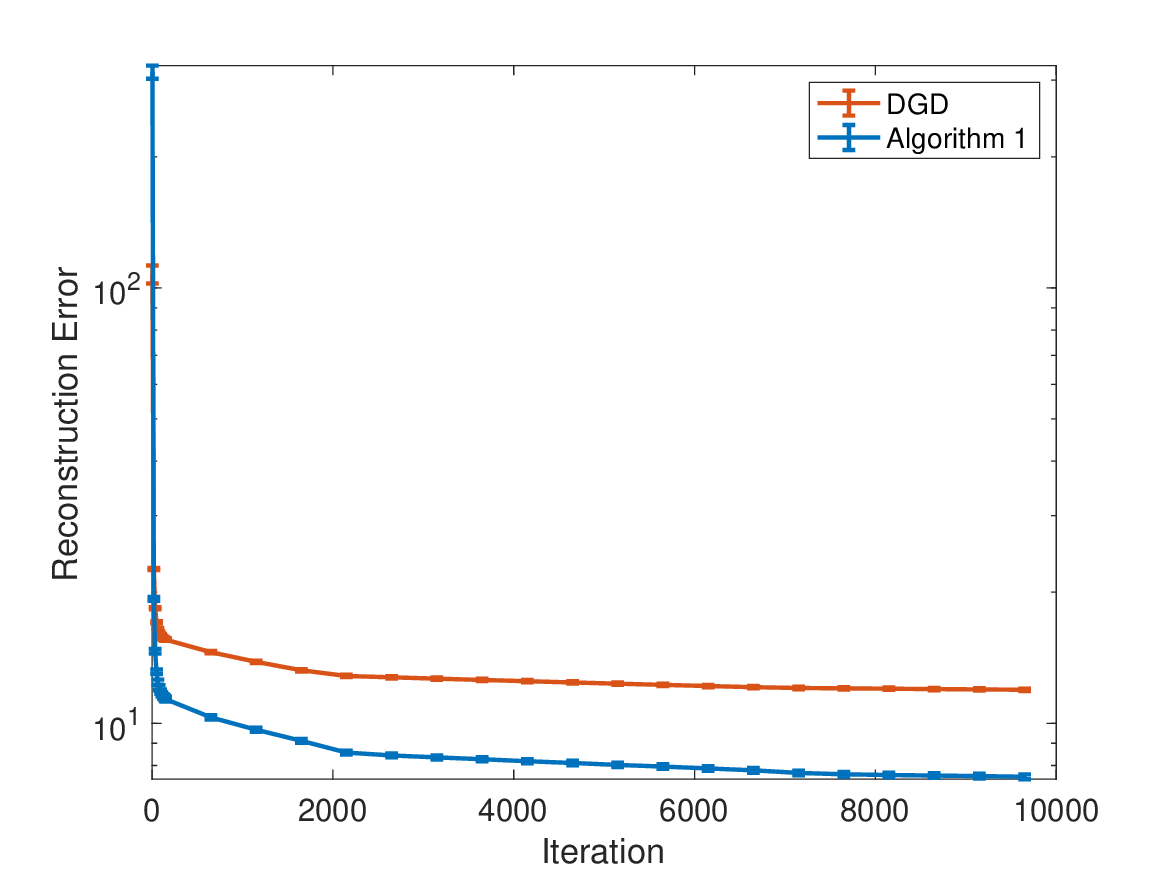}
\caption{Comparison of the objective function value between the proposed Algorithm \ref{alg:1} and the existing algorithm DGD in \cite{yuan2016convergence}. }
\label{Objpng}
\end{figure}

\begin{figure}
\centering
\includegraphics[width=.8\linewidth]{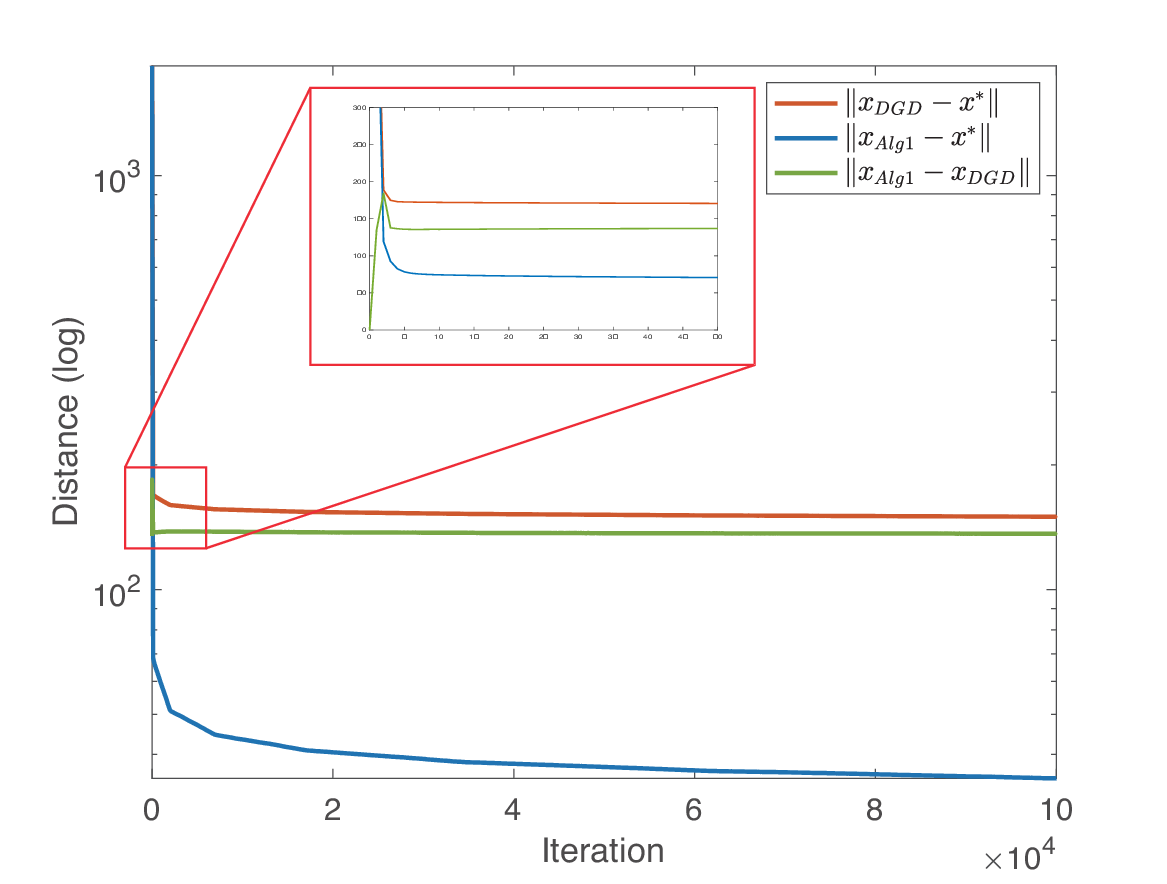}
\caption{Comparison of the distance between learned parameters and the actual optimal solution $x^\ast$ (obtained using centralized optimization). The learned parameters in our algorithm are represented as $ x_{Alg1} $, and the learned parameters in the existing algorithm DGD are represented by $x_{DGD}$. It can be seen that our algorithm does converge to a better solution than DGD. }
\label{distancepng}
\end{figure}

Fig.~\ref{Objpng} shows the evolution of the objective function values under our algorithm and the existing algorithm DGD in \cite{yuan2016convergence}, respectively. It is clear that our algorithm gives a much smaller cost value. To show that this is indeed due to different convergence properties between our algorithm and DGD, in Fig.~\ref{distancepng}, we plot the distance between learned parameters and the global optimal parameter, which is obtained using centralized optimization. It is clear that our algorithm indeed converges to a much better solution than DGD, likely due to its ability to evade saddle points.

\subsection*{Convolutional Neural Network}
\begin{figure}[tbhp]
\centering
\includegraphics[width=.8\linewidth]{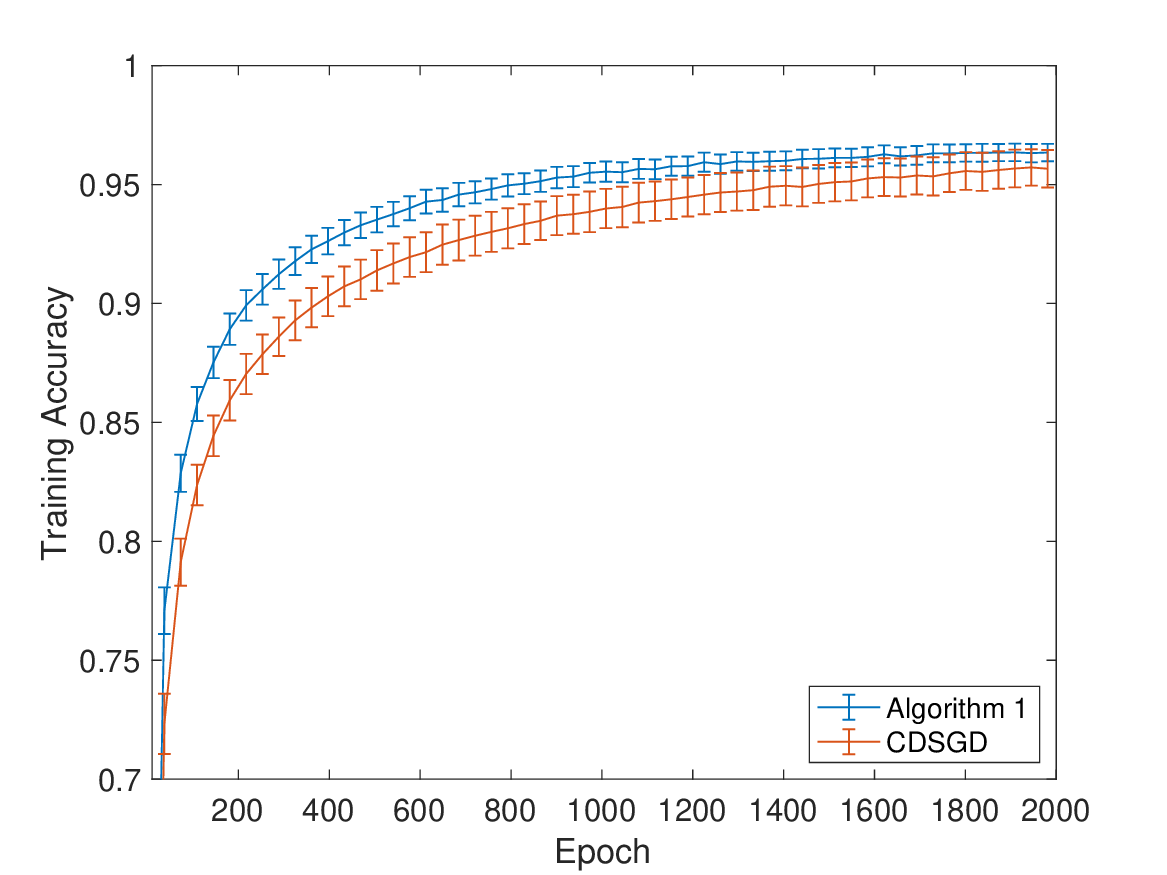}
\caption{Comparison of training accuracy between the proposed algorithm and a commonly used algorithm CDSGD from \cite{jiang2017collaborative}. }
\label{Fig:Train}
\end{figure}
\begin{figure}[tbhp]
\centering
\includegraphics[width=.8\linewidth]{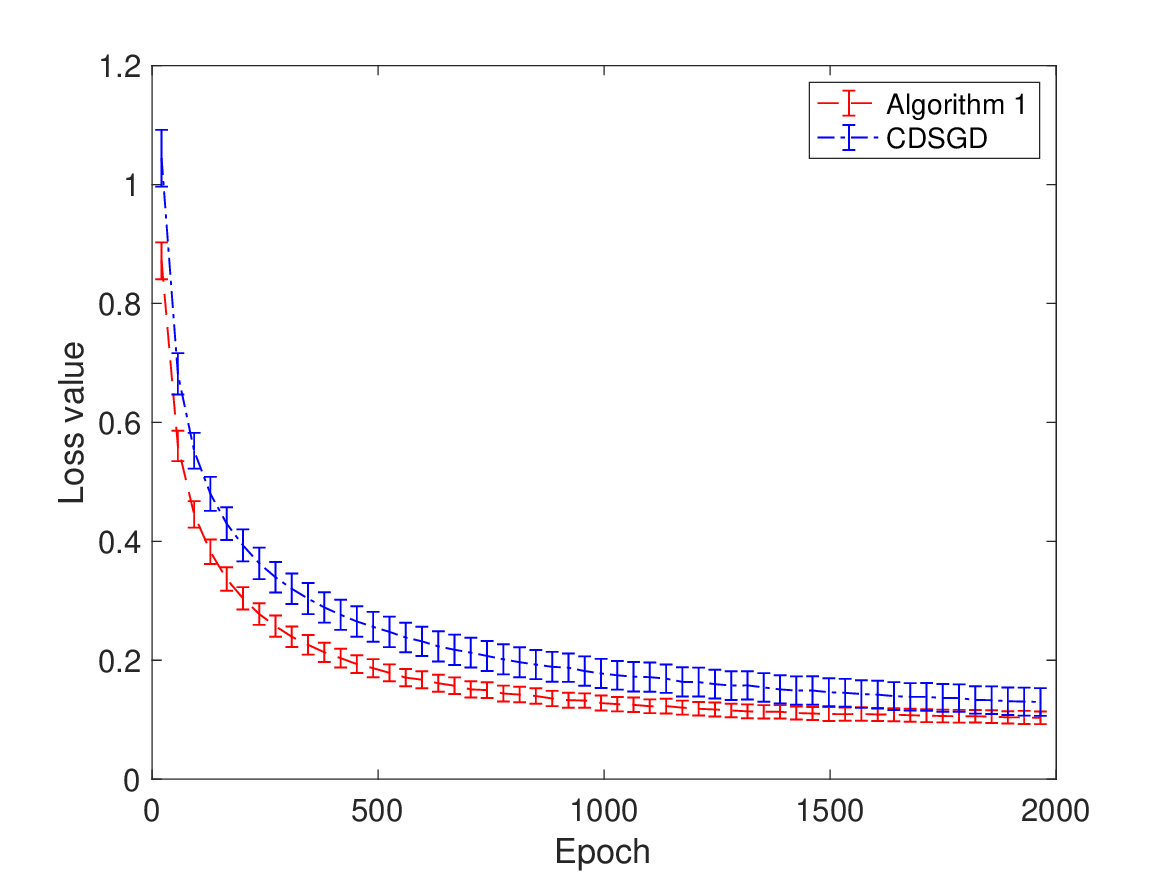}
\caption{Comparison of loss function value between the proposed algorithm and a commonly used algorithm CDSGD from \cite{jiang2017collaborative}. }
\label{Fig:Loss}
\end{figure}
For this experiment, we consider the training of a convolutional neural network (CNN) for the classification of the CIFAR-10 dataset, which contains 50,000 training images across 10 different classes. We evenly spread the CIFAR-10 dataset to the five agents, and set the batch size as 32. Our baseline CNN architecture is a deep network ResNet-18, the training of which is a highly nonconvex problem characterized by the presence of many saddle points \cite{dauphin2014identifying}. In the experiments, we train the CNN using both the proposed Algorithm \ref{alg:1} and the decentralized optimization algorithm CDSGD proposed in \cite{jiang2017collaborative}. In order to ensure fairness in comparison, both algorithms use the same set of learning rates ($\alpha= 0.62$, $\beta = 0.94$, $c_1 = 0.5$, $c_2 = 0.3$). The quantization interval $\ell$ is set such that all quantized outputs are representable using a binary string of $10$ bits.

The evolution of the training accuracies and loss-function values averaged over 10 runs are illustrated in Fig.~\ref{Fig:Train} and Fig.~\ref{Fig:Loss}, respectively.
It is evident that Algorithm \ref{alg:1} achieves lower loss function values more rapidly compared to CDSGD. This difference indicates that controlled quantization effects in our algorithm can aid in evading saddle points and discovering better function values.

\subsection*{Tensor Decomposition}
In this experiment, we consider Tucker tensor decomposition on the neural dataset in \cite{TDdataset}.  For $N$ neurons over $K$ experimental trials, when each trial  has $T$ time samples, the recordings of firing activities can be represented as an $N \times T \times K$ array, which is also called a third-order tensor \cite{williams2018unsupervised}. Each element in this tensor, $x_{n,t,k}$, denotes the firing rate of neuron $n$ at time $t$ within trial $k$.   Tucker tensor decomposition decomposes a tensor into a core tensor multiplied by a matrix along each mode. Following \cite{williams2018unsupervised},  we consider the tensor decomposition problem  for a  tensor recording  $\mathscr{X} \in \mathbb{R}^{50 \times 500 \times 100}$ of neural firing activities:
\begin{equation}
\mathscr{X}\approx \mathcal{T} \times _1 A\times _2 B \times _3 C=\sum_{n=1}^N{\sum_{t=1}^T{\sum_{k=1}^K{t_{n,k,t}}a_n\circ b_t\circ c_k}},
\end{equation}
where  $\circ$ represents the vector outer product, $\times_i$ (with $i=\{1,\,2,\,3\}$) denotes the $i$-mode matrix product,
$\mathcal{T}\in \mathbb{R}^{5 \times 5 \times 5}$ is the core tensor, and  $A \in \mathbb{R}^{50 \times 5}$, $B \in \mathbb{R}^{500 \times 5}$ and $C \in \mathbb{R}^{100 \times 5}$ are the three factors for Tucker decomposition.
The goal of tensor decomposition is to minimize the normalized reconstruction error  $\mathcal{E} =\left( \parallel \mathscr{X}-\mathcal{T} \times _1A\times _2B\times _3C\parallel _{F}^{2} \right) /\parallel \mathscr{X}   \parallel _{F}^{2}$, where the subscript $F$ denotes the Frobenius norm. It is well known that the tensor decomposition problem is inherently
susceptible to the saddle point issue \cite{ge2015escaping}.

We implement both the DGD algorithm in \cite{yuan2016convergence} and our Algorithm  \ref{alg:1}  to solve the tensor decomposition problem. For the DGD algorithm, we use the largest constant stepsize that can still ensure convergence, and for our algorithm, we set the stepsize parameters as $\alpha= 0.61$, $\beta = 0.92$, $c_1 = 0.03$, and $c_2 = 0.3$. The quantization interval $\ell$ is set such that all quantized outputs are representable using a binary string of $6$ bits.

\begin{figure}[tbhp]
\centering
\includegraphics[width=.8\linewidth]{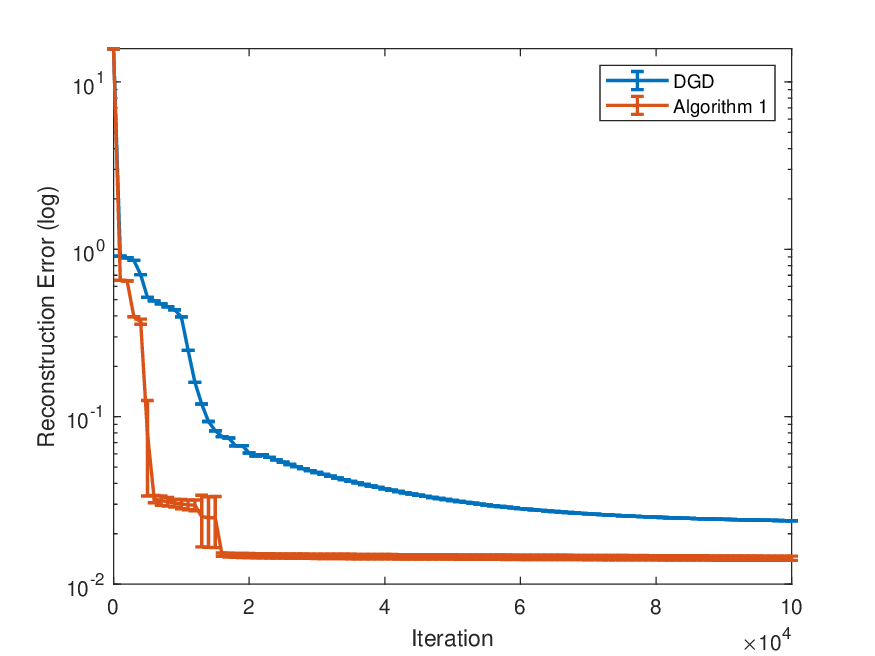}
\caption{Comparison of reconstruction error in tensor decomposition between the proposed Algorithm 1 and the existing algorithm DGD in \cite{yuan2016convergence}.}
\label{Fig:TD1}
\end{figure}
The evolution of the reconstruction error for the two algorithms under 50 runs is shown in Fig.~\ref{Fig:TD1}. It is clear that our algorithm finds better optimization solutions by effectively evading saddle points.

\subsection*{Robust Principal Component Analysis (PCA)}
In this experiment, we consider the problem of background subtraction in computer vision using robust PCA. Compared with the conventional PCA, robust PCA can provide a low-dimensional approximation that is more robust to outliers in data samples. For a given sequence of images (video), we employ robust PCA to separate moving objectives in the video from the static background. More specifically, for a given sequence of images represented as a data matrix
 $\boldsymbol{M} \in \mathbb{R}^{m \times n}$, we use robust PCA  to decompose $\boldsymbol{M}$ into a low-rank matrix $\boldsymbol{U}\boldsymbol{V}^\top$ (representing the background) and a sparse matrix $\boldsymbol{S}$ (representing moving objects), where $\boldsymbol{U} \in \mathbb{R}^{m \times r}$, $\boldsymbol{V} \in \mathbb{R}^{n \times r}$, $\boldsymbol{S} \in \mathbb{R}^{m \times n}$ and $ r \ll \min\{m, n\} $. Mathematically, the problem can be formulated as the following optimization problem \cite{ma2018efficient}:
\begin{equation}
\begin{aligned}\label{Equ:RPCA}
\min_{\boldsymbol{U}, \boldsymbol{V}} & f(\boldsymbol{U}, \boldsymbol{V}) + \mu_2 \|\boldsymbol{U}^\top\boldsymbol{U}-\boldsymbol{V}^\top\boldsymbol{V}\|^2_{F},\\
& f(\boldsymbol{U}, \boldsymbol{V}) = \min_{\boldsymbol{S} \in \mathcal{S}_{\bar{\alpha}}} \frac{1}{2} \| \boldsymbol{U}\boldsymbol{V}^\top + \boldsymbol{S} - \boldsymbol{M}\|^2_{F},
\end{aligned}
\end{equation}
where $\mu_2$ is a constant and  $\mathcal{S}_{\bar{\alpha}}$ represents the set of matrices with at most $\bar{\alpha}-$fraction of nonzero entries in every column and every row.

In our experiment, we use the ``WallFlower" datasets from Microsoft \cite{RPCAdataset}. We randomly assign 200 image frames with $56 \times 56$ pixels to each agent, resulting in the data matrix $\boldsymbol{M}_i$ of agent $i$ being of dimensions $m=9408$ and $n=200$.
 We set $\mu_2$ to 0.01, $\bar{\alpha}=0.2$, and $r=30$, and then solve [\ref{Equ:RPCA}] using the gradient descent based algorithm (``Fast RPCA") in \cite{yi2016fast}.  Fast RPCA employs a sorting-based estimator to generate an initial estimate $ \boldsymbol{S}_0$ and then it employs singular value decomposition  to generate the corresponding initial values of $\boldsymbol{U}_0$ and $\boldsymbol{V}_0$. Fast RPCA alternates between taking gradient steps for $\boldsymbol{U}$ and $\boldsymbol{V}$, and computing a sparse estimator to adjust $\boldsymbol{S}$.
In the experiment, we use  the best constant stepsize that we can find for   Fast RPCA (the largest stepsize that can still ensure convergence).  For our algorithm, we set the stepsize parameters  as $\alpha= 0.61$, $\beta = 0.92$, $c_1 = 0.003$, and $c_2 = 0.3$. The quantization interval $\ell$ is set such that all quantized outputs are representable using a binary string of $5$ bits.

\begin{figure}[tbhp]
\centering
\includegraphics[width=.8\linewidth]{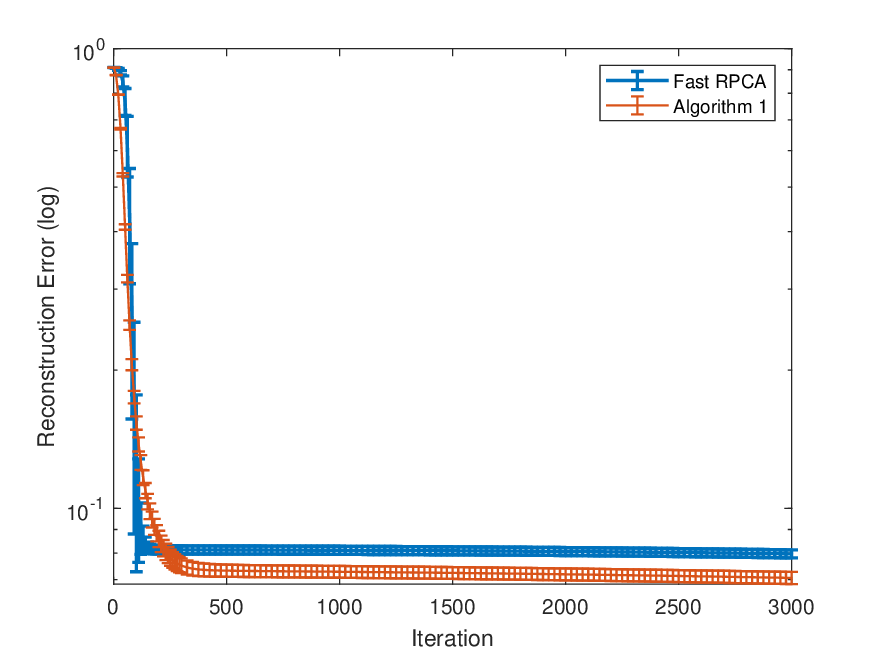}
\caption{Comparison of the reconstruction error in Robust PCA between the proposed Algorithm 1 and the existing algorithm Fast RPCA in \cite{yi2016fast}.}
\label{Fig:RPCA1}
\end{figure}
Fig.~\ref{Fig:RPCA1} shows the evolutions of the reconstruction error $\mathcal{E} =\sum_{i=1}^{N}   \| \boldsymbol{M}_i -\boldsymbol{U}_i\boldsymbol{V}_i^\top - \boldsymbol{S}_i \|  / \| \boldsymbol{M}_i \| _{F}^{2}$ under our algorithm and   Fast RPCA in \cite{yi2016fast}, respectively.
It is clear that our algorithm is capable of identifying superior solutions that yield a smaller reconstruction error. This implies that our algorithm can locate more favorable stationary points by effectively avoiding strict saddle points (it has been proven in \cite{ge2015escaping} that all saddle points in robust PCA are strict saddle points).

\section*{Discussions}

\subsection*{On Comparison with Other Stepsize Strategies}
To test if our stepsize strategy leads to a reduced convergence speed compared with existing counterparts which do not consider saddle-point avoidance, we also conduct experiments using the tensor decomposition problem to compare the convergence speed under our stepsize strategy, the constant stepsize strategy, a random stepsize strategy, and the conventional diminishing  stepsize strategy. For the
constant stepsize case, we use the largest constant stepsize that does not lead to divergence, and for the random stepsize strategy, we select the stepsize values in the ``hold" stages of our approach randomly from the reference functions $\frac{0.03}{1+0.6t^{0.61}}$ and $\frac{0.003}{1+0.6t^{0.92}}$. For the diminishing
stepsize case, we use the reference functions as the stepsizes,  which are commonly used in distributed optimization. The simulation results in Fig. S1 of the Supporting Information show that our algorithm can provide similar or even faster convergence speeds, and hence show that our approach does not trade convergence speed for saddle-point avoidance.

\subsection*{On Comparison with the Log-scale Quantization}
It is worth noting that recently \cite{9904851} and \cite{doostmohammadian2022fast} propose to use log-scale quantization in distributed optimization and prove that accurate convergence can be ensured when the objective functions are convex. However, the log-scale quantization scheme is not appropriate for the saddle-point avoidance problem in distributed nonconvex optimization. This is because to enable saddle-point avoidance, we have to keep the magnitude of quantization error large enough to perturb the optimization variable, no matter what the value of the optimization variable is (because we do not know where the saddle-point is). In fact, this is why we introduce the periodic switching between two sets of quantization levels in our quantization scheme (to avoid the possibility that a quantization input coincides with an endpoint of a quantization interval and results in a zero quantization error). However, the log-scale quantization scheme results in a quantization error that can be arbitrarily small when the quantization input is arbitrarily close to zero, meaning that the quantization-induced perturbation becomes negligible when the quantization input is close to zero, making it inappropriate for saddle-point avoidance. In fact, our experimental results using the binary classification problem in  Fig. S2 of the Supporting Information also confirm that the log-scale quantization scheme cannot provide comparable performance with our proposed quantization scheme.

\subsection*{On Applicability to High-Order Optimization Methods}
Given that additive noises have been proven effective in evading saddle points in second-order optimization algorithms as well (see, e.g., \cite{paternain2019newton}), our quantization effect based approach is well positioned to help saddle-point avoidance in second-order nonconvex optimization algorithms. To confirm this point,   we apply the quantization scheme to second-order Newton-method based distributed optimization  for the binary classification problem     (see details in the section ``Experimental Results Based on the Newton Method" on page 19 of the Supporting Information). The results in Fig. S3 in the Supporting Information confirm that our quantization scheme does significantly enhance the quality of the solution by evading saddle points compared with the case without  quantization effects. We plan to systematically investigate exploiting quantization effects in high-order    optimization algorithms to evade saddle points in future work.

\subsection*{On Relaxing the Smoothness Assumption} In the theoretical analysis, we assume that the objective functions are Lipschitz continuous. Given that ``generalized gradients" \cite{cortes2008discontinuous} have been proven effective to address non-smooth objective functions in convex optimization, it is tempting to investigate if the generalized gradient approach can be exploited to address nonconvex and non-smooth objective functions. Unfortunately, \cite{zhang2020complexity} proves that in general nonconvex and non-smooth optimization, for any $\epsilon\in[0,1)$, there is a more than $ 50\%$ probability that an $\epsilon$-first-order-stationary point (defined in the sense of the generalized gradient, usually called Clarke stationary point)  can never be found by any finite-time algorithm.    In future work, we plan to explore if some subclasses of nonconvex and non-smooth objective functions can be addressed using the generalized gradient approach.
\section*{Conclusions}

Saddle-point avoidance is a fundamental problem in nonconvex optimization. Compared with the centralized optimization case, saddle-point avoidance in distributed optimization faces unique challenges due to the fact that individual agents can only access local gradients, which may be significantly different from the global gradient (which actually carries information about saddle points). We show that quantization effects, which are
unavoidable in any digital communications, can be exploited without additional cost to evade saddle points in distributed nonconvex optimization. More specifically, by judiciously co-designing the quantization scheme and the stepsize strategy, we propose an algorithm that can ensure saddle-point avoidance and convergence to second-order stationary points in distributed nonconvex optimization.  Given the widespread applications of distributed nonconvex
optimization in numerous engineered systems and deep learning, the results are expected to have broad ramifications in various fields involving nonconvex optimization.
Numerical experimental results using distributed optimization and learning applications on benchmark datasets confirm the effectiveness of the proposed algorithm.

\showmatmethods{} 

\acknow{The work was supported in part by the National Science Foundation
under Grants CCF-2106293, CCF-2215088, and CNS-2219487.}

\showacknow{} 

\bibsplit[15]

\bibliography{pnas-sample}

\end{document}